\documentclass[letterpaper,twocolumn]{autart} 

\usepackage{filecontents}
\usepackage{calc}
\usepackage{latexsym}
\usepackage{amsfonts}
\usepackage{mathrsfs}
\usepackage{algorithm}
\usepackage{color}
\usepackage{url}
\usepackage{amssymb,amsmath,amsfonts,layout}
\usepackage{makeidx}
\usepackage{blkarray,multirow}
\usepackage{cite}
\usepackage{enumerate}
\usepackage{algpseudocode}
\usepackage{braket}
\usepackage{tikz}
\usepackage{lipsum}
\usepackage{caption}
\usepackage{subcaption}

\DeclareMathOperator{\R}{\mathbb{R}}
\DeclareMathOperator{\Hc}{\mathcal{H}}

\newcommand{\real}{\ensuremath{\mathbb{R}}}

\newcommand{\norm}[1]{\left\lvert\left\lvert#1\right\rvert\right\rvert}

\newtheorem{assumption}{Assumption}
\newtheorem{theorem}{Theorem}[section]

\newtheorem{remark}[theorem]{Remark}
\newtheorem{definition}[theorem]{Definition}

\newtheorem{problem}{Problem}

\newcommand{\longthmtitle}[1]{\mbox{}{\bf \textit{(#1).}}}
\allowdisplaybreaks

\newcommand{\rev}[1]{{\color{blue} #1}}
\renewcommand{\rev}[1]{#1}

\begin{document}
	\begin{frontmatter}
		
		\title{Online Data-Enabled Predictive Control} 
		
		\thanks{This work was authored in part by NREL, operated by Alliance for Sustainable Energy, LLC, for the U.S. Department of Energy (DOE) under Contract No. DE-AC36-08GO28308. Funding provided by DOE Office of Electricity, Advanced Grid Modeling Program, through agreement NO. 33652. }
		\thanks[footnoteinfo]{The first two authors contributed equally.}
		
		\author[NREL]{Stefanos Baros}\ead{stefanos.baros@nrel.gov},
		\author[NREL]{Chin-Yao Chang}\ead{chinyao.chang@nrel.gov},    
		\author[UCB]{Gabriel E. Col\'on-Reyes}\ead{gecolonr@berkeley.edu},   
		\author[NREL]{Andrey Bernstein}\ead{andrey.bernstein@nrel.gov}  
		
		\address[NREL]{National Renewable Energy Laboratory (NREL), Golden, CO 80401, USA}  
		\address[UCB]{University of California, Berkeley, CA 94720, USA} 
		
		\begin{abstract}
			We develop an online data-enabled predictive (ODeePC) control method for optimal control of unknown systems, building on the recently proposed DeePC \cite{DeePC}. Our proposed ODeePC method leverages a primal-dual algorithm with real-time measurement feedback to iteratively compute the corresponding real-time optimal control policy as system conditions change. 
			\rev{The proposed ODeePC conceptual-wise resembles standard adaptive system identification and model predictive control (MPC), but it provides a new alternative for the standard methods. ODeePC is enabled by computationally efficient methods that exploit the special structure of the Hankel matrices in the context of DeePC with Fast Fourier Transform (FFT) and primal-dual algorithm. } 
			We provide theoretical guarantees regarding the asymptotic behavior of ODeePC, and we demonstrate its performance through \rev{numerical examples}.  
		\end{abstract}
	\end{frontmatter}
	
	\section{Introduction}
	In the context of smart critical infrastructures and complex dynamic systems,  applications of optimal  control abound \cite{Baros}. Traditional optimal control of these systems relies on accurate system models, with model predictive control (MPC) \cite{Dorfler_4,Dorfler_5,Dorfler_6,Dorfler_7,Dorfler_8} being a classic approach. Many complex physical systems however, such as large power systems, manifest complex and hard-to-model uncertain dynamics that complicate their study and analysis. Therefore, data-based modeling and control techniques have become increasingly popular in the context of such systems due to data abundance \cite{Dorfler_1, Dorfler_2}.
	Such approaches offer an attractive alternative to classic optimal control as they are independent of any analytical system models.
	System identification algorithms are exploited in many control approaches for producing approximate models  when an analytical system description  might not be available or might be too difficult to obtain \cite{Dorfler_22, Dorfler_23, Dorfler_24}. Other data-based control approaches exploit different types of learning algorithms at various stages of the control design process to generate information about the model. Representative publications of this line of work are \cite{Dorfler_16, Dorfler_17,Dorfler_18,Dorfler_19,Dorfler_20,Dorfler_21}, \cite{Dorfler_28}, and \cite{Dorfler_29}; the most prominent approaches are based on reinforcement learning and MPC using Gaussian processes. All these approaches rely on    state-space model representation of the underlying dynamic system.

	Data-driven optimal control approaches can be useful even when a certain (simplified) model of the system exists. In such a case,  data-driven approaches can help refine the model, by  capturing  exogenous disturbances and other model details not accounted for in the existing model \cite{Dorfler_28}. \rev{The idea of using series of input/output data to characterize system behavior was first explored in 
		\cite{willems1986time_1,willems1986time_2,willems1987time_3}. Since then, the concept had been matured with behavioral system fundamental theory  being established~\cite{willems2005note,markovsky2006exact}. The concept was also explored in various 
		model-free control schemes~\cite{favoreel1999model,kadali2003data,markovsky2008data}. Recently, \cite{DeePC} revived the behavioral system approach by framing so-call data-enabled predictive control algorithm (DeePC). It basically carries the spirit of totally data-based approach which uses the control inputs and plant outputs of a dynamic system to learn the system's behavior and compute a predictive control policy directly.
	}
	Despite its accuracy and efficacy for small systems (systems with few states), however, this control approach, lacks scalability as its implementation to large systems (i.e., systems with many states) can be computationally burdensome. This is partially because of the large dataset required to accurately represent the underlying system.
	
	\rev{
		\textbf{\textit{Contributions}}. In this paper, we develop an \emph{online} data-enabled predictive controller (ODeePC) built on DeePC. The objective of ODeePC is to utilize the online data measurements to capture real-time variation of the controlled system, and update the control accordingly. Conceptually, ODeePC  is similar to standard methods of MPC with online parameter/model identification. We view ODeePC as a valuable alternative for predictive control, especially for black box systems that are challenging for standard MPC-based approaches. The capability of ODeePC is enabled by the following elements.
		\begin{itemize}
			\item To enable online implementation, we leverage the online time-varying optimization methods with measurement feedback \cite{AndreyOnlineOpt,bernstein2019real}. Under this setup, the measurement is not directly used to update the decision variables; instead, it is used to update the behavioural system model and, in turn, adapt the control implicitly. We also use the well-established method for proof of convergence for online optimization~\cite{AndreyOnlineOpt,bernstein2019real} for the convergence of ODeePC.
			\vspace{2mm}
			\item ODeePC requires computing products of large-scale matrices and vectors. The computational burden may hinder the online implementation of the algorithm. We devise a computationally efficient algorithm for a fast computation of the product of a non-square block Hankel matrix with a vector. The algorithm exploits properties of circulant matrices and fast Fourier transform (FFT) to carry out the assigned computations efficiently. We also derive the complexity of this algorithm and prove that it precisely computes the desired product.
		\end{itemize}
	}
	
	The rest of this paper is structured as follows. In Section \ref{sec:probform}, we present the notation we use throughout the paper and some preliminaries. In Section \ref{sec:MPC_DeePC}, we review the classic MPC and the recently proposed DeePC \cite{DeePC}. Section \ref{sec:ODeePC} includes the main results surrounding the development of the proposed ODeePC algorithm and the analysis of its performance. We validate the performance of ODeePC through numerical examples in Section \ref{sec:ODeePCVal}. Finally, in Section \ref{sec:conc}, we conclude the paper with some remarks.
	
	\section{Notation and Preliminaries}
	\label{sec:probform}
	\subsection{Notation}
	Let $\R$ denotes the set of real numbers; $\mathbb{Z}_{\geq0}$ and $\mathbb{Z}_{>0}$ respectively, denote the set of nonnegative integers. Given a matrix $A$, $A^\top$ denotes its transpose, and ${A\succ(\prec) 0}$ denotes that $A$ is positive (negative) definite. The matrix $I_n \in \R^{n\times n}$ is the $n \times n$ identity matrix. For $x\in\real^n$, $\norm{x}_2$ denotes its Euclidean norm, $\text{diag}\{x\}$ is a diagonal matrix with the elements of $x$ on the main diagonal. Further, given a set $\mathcal{X}\subset \real^n$, $\text{Proj}_{\mathcal{X}}\{x\}$ denotes the projection of $x$ onto $\mathcal{X}$. Given a $Q\succ 0$, we define $\norm{v}_Q^2 = v^\top Q v$. We use subscript $t$ to denote a vector value at time $t$; and $v_{i,t}$ to denote the $i^{th}$ entry of the vector $v_t$ at time $t$.  
	We denote the smallest eigenvalue of $Q\in\real^{n\times n}$ by $\lambda_{\min}(Q)$. 
	A shift matrix is defined as $S_{m,n} = \begin{bmatrix} 0 & I_{m-1} \\ 0 & 0 \end{bmatrix} \otimes I_n$ and $S_{m,n}\in\real^{mn\times mn}$. For any matrix $A \in \real^{mn \times N}$, $S_{m,n} A$ results in shifting the elements of $A$ upward by $n$ positions. 
	
	\subsection{Preliminaries}
	Let us start by assuming that the system we seek to control can be represented by a discrete-time, linear, time-invariant (LTI), state-space model given by
	\begin{equation}
	\label{eq:dss}
	\begin{cases}
	x_{t+1} = Ax_t + Bu_t \\    
	y_t = Cx_t + Du_t. 
	\end{cases}
	\end{equation}
	In this representation, $t \in \mathbb{Z}_{>0}$ is the discrete time index; 
	$A \in \mathbb{R}^{n \times n}$ is the state matrix; $ B \in \mathbb{R}^{n \times m}$ is the input matrix; $ C \in \mathbb{R}^{p \times n}$ is the output matrix;  $ D \in \mathbb{R}^{p \times m}$ is the feed-forward matrix;  $x_t\in\mathbb{R}^n$ is the state vector; $u_t\in\mathbb{R}^m$ is the control input vector; and $y_t\in\mathbb{R}^p$ is the output vector. Also, let $r_t\in\mathbb{R}^p$ denote the output \emph{reference} vector representing the desired value of the output at time $t$.
	
	Let $u:=(u_1^\top,...,u_T^\top)^\top \in \real^{mT}$ denote the vector of the control inputs during $t = 1, \ldots, T$, for some $T \in \mathbb{Z}_{>0}$. The \emph{block Hankel matrix}, with each column being a set of consecutive data points of length $L\le T$, is given by 
	\begin{equation}
	\label{eq:Hankel}
	\mathscr{H}_L(u) := \begin{bmatrix} u_1 && u_2 && \dots && u_{T-L+1} \\ u_2 && u_3 && \dots && u_{T-L+2} \\ \vdots && \vdots && \ddots && \vdots \\ u_L && u_{L+1} && \dots && u_T \end{bmatrix},
	\end{equation}
	with $\mathscr{H}_L(u) \in \real^{mL \times (T - L + 1)}$. A similar matrix can be constructed for $y$. It is useful at this point to introduce the following definition.
	\begin{definition} \label{def:pers}
		Let $L,T \in \mathbb{Z}_{>0}$ and $T \geq L$. Then, the signal $u:=(u_1^\top,...,u_T^\top)^\top$ is \emph{persistently exciting} of order $L$ if $\mathscr{H}_L(u)$, as defined in \eqref{eq:Hankel}, is of full row rank.
	\end{definition}
	Roughly speaking, a persistently exciting input is a rich enough input able to excite the system so that it generates an output that is representative of its behavior. 
	
	Moving forward, we assume that $u \in \mathbb{R}^{mT}$ is a persistently exciting input data set of order $T_{\text{tot}}$. With rich enough data to describe system~\eqref{eq:dss}, we can construct the following block Hankel matrices associated with the input and output data, respectively 
	\begin{align}
	\label{eq:Hankel_matrices}
	U:= \mathscr{H}_{T_{\text{tot}}}(u), \;\;Y:= \mathscr{H}_{T_{\text{tot}}}(y),    
	\end{align}
	where $U\in\R^{mT_{\text{tot}}\times \kappa}$,  $Y\in\R^{pT_{\text{tot}}\times \kappa}$, and $\kappa := T - T_{\text{tot}}+1$. Behavioral system theory \cite{behaviorsystemstheory} states that if system~\eqref{eq:dss} is controllable\footnote{The \textit{controllability} for LTI in the context of behavioural system theory is weaker than the one for control system theory (behavioural system theory is on the control of output $y$ and control system theory is on the state $x$), cf.~\cite[Definition 4.3]{DeePC}. We do not distinguish them here for simplicity.}, then the input-output pair $(u_{\text{tot}},y_{\text{tot}})$, where $u_{\text{tot}}\in\mathbb{R}^{mT_{\text{tot}}}$ and $y_{\text{tot}}\in\mathbb{R}^{pT_{\text{tot}}}$, is a trajectory of \eqref{eq:dss} if and only if there exists a $g \in \mathbb{R}^{\kappa}$, such that
	\begin{equation}
	\label{eq:bst}
	\begin{bmatrix} U \\ Y \end{bmatrix}g =
	\begin{bmatrix} u_{\text{tot}} \\ y_{\text{tot}} \end{bmatrix}.  
	\end{equation}
	Equation~\eqref{eq:bst} indicates that any possible trajectory should be a linear combination of $\kappa$ number of trajectories (columns) embedded in $[U^\top, Y^\top]^\top$. We can view~\eqref{eq:bst} as an alternative model of~\eqref{eq:dss}. The main difference is that~\eqref{eq:bst} is entirely constructed by data as opposed to state evolution in~\eqref{eq:dss}.
	\section{Overview of Classic MPC and DeePC}
	\label{sec:MPC_DeePC}
	In this section, we provide an overview of the classic MPC and the recently proposed DeePC \cite{DeePC}. Note that the classic MPC uses a precise model~\eqref{eq:dss}, whereas DeePC uses the recorded control inputs and plant outputs of the underlying system to capture its behavior and compute a predictive control policy. 
	\subsection{MPC}
	MPC is a receding time horizon control algorithm that computes an optimal control input $u_t$ based on a prediction of the system's future trajectory subject to the system's dynamics. We consider a classic MPC setting with prediction horizon $N \in \mathbb{Z}_{>0}$ formulated as
	\begin{align}\label{eq:MPC}
	& \underset{\substack{x, u_k \in \hat{\mathcal{U}}, y_k \in \hat{\mathcal{Y}}, \\ \forall k = 0,\cdots N-1}}{\text{minimize}}
	\hspace{2mm}\sum_{k=0}^{N-1} f(u_k,y_k), \\ \nonumber
	& \text{s. t. }
	\hspace{5mm} x_{k+1} = Ax_k + Bu_k,\quad \forall k \in \{0, \cdots, N-1 \}, \\ \nonumber
	&\hspace{12mm} y_k = Cx_k + Du_k, \quad \forall k \in \{0, \cdots, N-1 \} , \\ \nonumber
	& \hspace{12mm} x_0=\hat{x}_t ,
	\end{align}
	where $x=(x_0^\top,...,x_N^\top)^\top$; $\hat{\mathcal{U}} \subseteq \mathbb{R}^{m}$ and $\hat{\mathcal{Y}} \subseteq \mathbb{R}^{p}$ are, respectively, the convex and bounded constraint sets for $u_k$ and $y_k$ for all $k$; $f:\hat{\mathcal{U}}\times\hat{\mathcal{Y}} \mapsto \R$ is a convex cost function; and $\hat{x}_t$ serves as the initial state of the system. 
	We let  $t$ denote the current time, and  $k$ index the time instances of the look-ahead horizon window. Algorithm  \ref{MPC_algor} \cite{DeePC} can be used to solve the MPC Problem \eqref{eq:MPC}.
	\begin{algorithm}
		\underline{Problem data:} matrices $A,B,C,D$, current state $\hat{x}_t$, feasible input and output sets $\hat{\mathcal{U}}$ and $\hat{\mathcal{Y}}$, objective function $f$.
		\begin{enumerate}
			\item Solve \eqref{eq:MPC} for $u^\star=(u_0^{\star^\top},...,u_{N-1}^{\star^\top})^\top$.
			\item Apply inputs $(u_t^\top,...,u_{t+s}^\top)^\top=(u_0^{\star^\top},...,u_s^{\star^\top})^\top$ for some $s\leq N-1$.
			\item Set $t$ to $t+s$ and update $\hat{x}_t$.
			\item Repeat.
		\end{enumerate}
		\caption{MPC \cite{DeePC}}
		\label{MPC_algor}
	\end{algorithm}
	Usually, after the optimal control policy is computed, only the input $u_0^\star$ that corresponds to the first look-ahead window is implemented in the system (or $s=0$). In the rest of the paper, we assume $s=0$ without loss of generality.
	
	The  MPC algorithm has been proven to be effective in numerous applications, e.g.,  
	autonomous driving \cite{Dorfler_9} and flight control \cite{Dorfler_10}, where the goal is primarily trajectory tracking.  Despite that, the requirement for an accurate model description still restricts the application domain as systems whose dynamics are hard to model or unknown and cannot be considered. To this end,  DeePC, which leverages measured system data instead of an accurate system model to capture the system's behavior, overcomes these limitations and has been shown to work well for small systems. We review the DeePC control approach \cite{DeePC} next.
	
	\subsection{DeePC}
	DeePC relies on the past input/output data to construct the model shown in~\eqref{eq:bst}.
	Let $u_{\text{ini}}\in\R^{mT_\text{ini}}$ and $y_{\text{ini}}\in\R^{pT_\text{ini}}$ denote a given initial trajectory of the system of length $T_\text{ini}\in\mathbb{Z}_{>0}$ over time interval $[t-T_\text{ini},\cdots,t-1]$.  Any trajectory $u=(u_0^\top,...,u_{N-1}^\top)^\top$ and  $y=(y_0^\top,...,y_{N-1}^\top)^\top$ over the time interval $[t, t+N-1]$ should satisfy
	\begin{align}\label{eq:DeePC_model}
	\begin{bmatrix} U \\ Y \end{bmatrix}g = \begin{bmatrix} U_p \\ U_f \\ Y_p \\ Y_f \end{bmatrix}g =
	\begin{bmatrix} u_{\text{ini}} \\ u \\ y_{\text{ini}} \\ y \end{bmatrix},
	\end{align}
	for some $g\in\R^{\kappa}$. Note that we have $T_{\text{tot}} = T_{\text{ini}} + N$; and we have partitioned $U\in\R^{mT_{\text{tot}}\times \kappa}$ 
	into $U_p$ and $U_f$, with $U_p \in \R^{mT_{\text{ini}}\times \kappa}$ and $U_f \in \mathbb{R}^{mN\times \kappa}$.
	Similarly, $Y\in\R^{pT_{\text{tot}}\times \kappa}$ is partitioned into $Y_p \in \mathbb{R}^{pT_{\text{ini}}\times \kappa}$ and $Y_f \in \mathbb{R}^{pN\times \kappa}$. The length of the initial trajectory, $T_\text{ini}$, should be selected large enough to ensure unique $y$ for any given $u$ cf.~\cite[Lemma 1]{behaviorsystemstheory}. \rev{In fact,  $T_\text{ini}$ should be  the number of data points required to ensure the system behavior being ``observed'' (similar to the observability condition in the context of classical control theory); see \cite{DeePC} for a  formal lower bound on $T_\text{ini}$.} We assume that the total number of measured input/output pairs, $T$, is large enough to construct a persistently exciting $U$ of order $T_{\text{tot}}$; cf.~Definition \ref{def:pers}. In the following, we use the shorthand notation $\mathcal{H} = [U_p^\top,  U_f^\top, Y_p^\top, Y_f^\top]^\top$.
	
	With the elements in place, we state the DeePC formulation of~\eqref{eq:MPC} as:
	\begin{align}\label{eq:DeePC}
	& \underset{g \in \mathbb{R}^\kappa, u \in \mathcal{U}, y \in \mathcal{Y}}{\text{minimize}}\sum_{k=0}^{N-1}f(u_k,y_k) \text{ s.t. \eqref{eq:DeePC_model} holds},
	\end{align}
	where $\mathcal{U}$ and $\mathcal{Y}$ are, respectively, the Cartesian products of $N$ number of $\hat{\mathcal{U}}$ and $\hat{\mathcal{Y}}$. 
	In \cite{DeePC}, the following algorithm is proposed for solving  \eqref{eq:DeePC}.
	\begin{algorithm} 
		\caption{DeePC \cite{DeePC}}
		\underline{Problem data:} past input data $u_{\text{ini}},y_{\text{ini}}$, $\mathcal{H} = [U_p^\top,  U_f^\top, Y_p^\top, Y_f^\top]^\top$, feasible input and output sets $\mathcal{U}$ and $\mathcal{Y}$, objective function $f$.
		\begin{enumerate}
			\item  Solve \eqref{eq:DeePC} for $g^\star$ and compute $u^\star=U_f g^\star$.
			\item Apply control inputs $(u_t^\top,...,u_{t+s}^\top)^\top=(u_0^{\star^\top},...,u_s^{\star^\top})^\top$ for some $s\leq N-1$.
			\item Set $t$ to $t+s$ and update past input/output data $u_{\text{ini}}$ and $y_{\text{ini}}$ with the $T_{\text{ini}}$ most recent data obtained from measurements.
			\item Repeat.
		\end{enumerate}\label{alg:deepc}
	\end{algorithm}
	
	We note that in Algorithm \ref{alg:deepc}, the matrix $\mathcal{H}$ is not updated over time as the data is recorded offline and used online to solve the receding horizon predictive control problem. Only the RHS elements $h = [u_{\text{ini}}^\top, u^\top, y_{\text{ini}}^\top, y^\top]^\top$ are updated as $t$ evolves. In cases where we are dealing with time-varying or nonlinear systems, not updating the matrix $\mathcal{H}$ over time might result in bad system representation and subsequently invalid computed control policies. We next present our proposed Online Data-enabled Predictive Control (ODeePC) which tackles these challenges.
	
	
	\section{Online Data-enabled Predictive Control}
	\label{sec:ODeePC}
	Below are the main aspects of the proposed algorithm; the details are given in the ensuing sections.
	\begin{itemize}
		\item It exploits available real-time data obtained from system measurements to dynamically update both the matrix $\mathcal{H}$ and the vector $h$.
		\item It uses a primal-dual gradient descent algorithm to iteratively compute the optimal control policy, thus allowing the intermediate control inputs to be implemented in the system in real-time. Further, it allows real-time measured information about the system's state to take part in the algorithm and affect the computed optimal control policy.
		\item It exploits a FFT-based algorithm to efficiently compute the products of block Hankel matrices with vectors.
	\end{itemize}
	These aspects of ODeePC will be explained in detail.  
	
	\subsection{Primal-dual Algorithm for the Regularized DeePC Problem}
	We first design a primal-dual algorithm to iteratively solve the Lagrangian formulation of~\eqref{eq:DeePC}. Then, we appropriately modify the designed algorithm  to arrive at the ODeePC algorithm's iterative update rule. 

	We start with considering the following min-max optimization problem associated with \eqref{eq:DeePC}:
	\begin{equation}
	\label{eq:min-max}
	\begin{aligned}
	& \underset{\nu}{\text{maximize }}\bigg(\underset{g \in \mathbb{R}^\kappa, u \in \mathcal{U}, y \in \mathcal{Y}}{\text{minimize }}\mathcal{L}(u,y,g,\nu)\bigg),
	\end{aligned}
	\end{equation}
	where $\nu\in\real^{N_{\nu}}$, $N_{\nu} = (m+p)(T_{\text{ini}}+N)$, and $\mathcal{L}$ is the \textit{regularized Lagrangian}, given as:
	\begin{equation}
	\label{eq:Lag}
	\begin{aligned}
	\mathcal{L}(u,y,g,\nu) &= \sum_{k=0}^{N-1}f(u_k,y_k) + \frac{\epsilon_g}{2}\norm{g}_{2}^2 \\ \nonumber
	&\hspace{10mm}+ \nu^\top(\mathcal{H}g - h) - \frac{\epsilon_{\nu}}{2}\norm{\nu}_{2}^2,
	\end{aligned}
	\end{equation}
	where $\epsilon_g > 0$ and $\epsilon_{\nu} > 0$ are  (small) constants. The regularization terms improve the convergence rate of the gradient-based method at the expense of converging to a point that is close to but not exactly the real optimal point of \eqref{eq:DeePC}.
	This particular regularization is widely used for solving convex optimization problems  \cite{AndreyOnlineOpt}, \cite{bernstein2019real}. Using the primal-dual gradient descent algorithm to solve \eqref{eq:min-max}, we arrive at the problem's saddle-flow dynamics. These, are given by
	\begin{subequations}
		\label{eq:sf_solved}
		\begin{align}
		u^{\tau+1} &= \text{Proj}_{\mathcal{U}}\{u^{\tau} -\alpha(\nabla_u{\tilde{f}}|_{(u^{\tau},y^{\tau})} - \nu_{u}^{\tau})\}, \label{eq:sf_solveda} \\
		y^{\tau+1} &= \text{Proj}_{\mathcal{Y}}\{y^{\tau} -\alpha\big(\nabla_y{\tilde{f}}|_{(u^{\tau},y^{\tau})} - \nu_{y}^{\tau}\big)\}, \label{eq:sf_solvedb} \\
		\label{eq:sf_solvedc}
		g^{\tau+1} &= g^{\tau} -\alpha(\mathcal{H}^\top\nu^{\tau} + \epsilon_g g^{\tau}), \\
		\label{eq:sf_solvedd}
		\nu^{\tau+1} &= \nu^{\tau} +\alpha(\mathcal{H}g^{\tau} - h -\epsilon_{\nu}\nu^{\tau}),
		\end{align}
	\end{subequations}
	where $\tau$ is the iteration number, $\alpha\in\real_+$ is the step size, and $\tilde{f} (u, y) :=\sum_{k=0}^{N-1}f(u_k,y_k)$.  
	Further, we define $\nu_{u}^{\tau}$ and $\nu_{y}^{\tau}$ as the elements of $\nu^{\tau}$ associated with the inputs and outputs, respectively. The saddle-flow dynamics~\eqref{eq:sf_solved} solve the \textit{static} optimization problem \eqref{eq:min-max}.
	
	Relatively large systems with numerous states would give rise to a large matrix $\mathcal{H}$. This would inevitably render the algorithm computationally very expensive. To see this, consider that the product of a general matrix $\mathcal{H}\in \mathbb{R}^{n\times m}$ with a $m-$element vector  (at every iteration as imposed by \eqref{eq:sf_solvedc} and \eqref{eq:sf_solvedd}), can be computed by carrying out  $n \times m$ multiplications and $(n \times m-n)$ additions. One can realize that exact and real-time computation of this product and thus the implementation of the algorithm would be quite challenging. This would be even more challenging in an online setting where,  $u_{\text{ini}}$, $y_{\text{ini}}$ and $\mathcal{H}$ would be frequently updated. Motivated by this, we carefully design ODeePC so that it is computationally efficient and practically implementable. We accomplish this by exploiting a computationally efficient algorithm that leverages FFT to compute the block Hankel matrix-vector multiplication quickly. The details will be provided in the sequel.
	
	\subsection{ODeePC}
	We now present the underlying time-varying optimization problem associated with our proposed ODeePC. 
	We consider, the vector $h$ in optimization~\eqref{eq:min-max} to be frequently updated with the latest input-output pair, $u_{\text{ini}}$ and $y_{\text{ini}}$. In addition, to allow our algorithm to cope with time-varying systems, we also update $\Hc$ at every iteration for better system characterization. The following formulation captures the time-varying properties described. 
	\begin{equation}
	\label{eq:ODeePC}
	\begin{aligned}
	& \underset{\nu}{\text{maximize }}\bigg(\underset{g \in \mathbb{R}^\kappa, u \in \mathcal{U}, y \in \mathcal{Y}}{\text{minimize }}\mathcal{L}^t(u,y,g,\nu)\bigg),
	\end{aligned}
	\end{equation}
	where $t$ captures the time instances when the optimization problem is updated. The time varying Lagrangian is the same as~\eqref{eq:min-max} except that $\Hc$ and $h$ are time dependent, given as
	\begin{align*}
	\Hc^t = \begin{bmatrix} U_p^t \\ U_f^t \\ Y_p^t \\ Y_f^t \end{bmatrix}, \quad h^t = \begin{bmatrix} u_{\text{ini}}^t \\ u \\ y_{\text{ini}}^t \\ y \end{bmatrix}.
	\end{align*}
	\begin{remark}\longthmtitle{ODeePC for linear time-varying (LTV) systems}
		{\rm For a LTI system, updating $\Hc$ online using the measured data is not very meaningful as the original $\Hc$ already captures all the properties that characterize the system; however, in order to allow ODeePC to deal with LTV systems or nonlinear dynamical systems, frequent online updating of the matrix $\Hc$ is necessary.}
	\end{remark}
	
	Let the variable $\tau$ track the algorithm iteration, and assume that the control is implemented on the system every time $N_I$ iterations have been completed. The iteration index $\tau$ and the system update instances $t$ are related as follows. When $\tau$ coincides with the system update instant $t$, then $\tau+N_I$ would coincide with the system update instant $t+1$. In our proposed ODeePC system, all the inner-loop iterations in the interval $(t,t+1)$ are carried out using the update rules~\eqref{eq:sf_solved}, with $\Hc$ being fixed and given by $\Hc^t = [{U_p^t}^\top \; {U_f^t}^\top \; {Y_p^t}^\top \; {Y_f^t}^\top]^\top$. In addition, the $u_{\text{ini}}$ and $y_{\text{ini}}$ elements of $h$ are, respectively, fixed at $u_{\text{ini}}^t$ and $y_{\text{ini}}^t$. Every $N_I$ iterations however, starting at instant $t+1$ (or $\tau+N_I$), these elements are updated, and ODeePC deploys the following update rules to compute the new input-output pairs and dual variables:
	\begin{subequations}
		\label{eq:sf_online}
		\begin{align}
		u^{\tau+1} &= \text{Proj}_{\mathcal{U}}\{\hat{u}^{\tau} -\alpha(\nabla_u{\tilde{f}}|_{(\hat{u}^{\tau},\hat{y}^{\tau})} - \hat{\nu}_{u}^{\tau})\}, \label{eq:sf_online-1} \\
		y^{\tau+1} &= \text{Proj}_{\mathcal{Y}}\{\hat{y}^{\tau} -\alpha\big(\nabla_y{\tilde{f}}|_{(\hat{u}^{\tau},\hat{y}^{\tau})} - \hat{\nu}_{y}^{\tau}\big)\}, \label{eq:sf_online-2} \\
		g^{\tau+1} &= g^{\tau} -\alpha(\Hc^{(t+1)^\top}\hat{\nu}^{\tau} + \epsilon_g g^{\tau}), \label{eq:sf_online-3}\\
		\nu^{\tau+1} &= \hat{\nu}^{\tau} +\alpha(\Hc^{t+1} g^{\tau} - h^{t+1} - \epsilon_{\nu} \hat{\nu}^{\tau}), \label{eq:sf_online-4}
		\end{align}
	\end{subequations}
	where $\hat{u}^{\tau} = S_{N,m}u^{\tau}$, $\hat{y}^{\tau} = S_{N,p}y^{\tau}$, $\hat{\nu}^{\tau} = [\hat{\nu}^{\tau^{\top}}_u\; \hat{\nu}^{\tau^{\top}}_y]^\top$, $\hat{\nu}^{\tau}_u = S_{T_{\text{tot}},m}\nu^{\tau}_u$,  $\hat{\nu}^{\tau}_y = S_{T_{\text{tot}},p}\nu^{\tau}_y$, and
	\rev{
		\begin{subequations}\label{eq:Update_H&h}
			\begin{align}
			& h^{t+1} = [u_{\text{ini}}^{t+1}\;\; u^{\tau+1}, \; y_{\text{ini}}^{t+1}\;\; y^{\tau+1} ],   \\ 
			& \Hc^{t+1} = \begin{bmatrix}S_{T_{\text{tot}},m} & 0 \\ 0 & S_{T_{\text{tot}},p} \end{bmatrix} \Hc^{t} + \Hc^{t}_{\text{add}}, \\ 
			&\Hc^{t}_{\text{add}} = \begin{bmatrix}
			0 \\ 
			[0\;\;I_m]U^{t}S_{\kappa,1}^{\top} +  \begin{bmatrix}
			0 & [I_m\;\; 0] u^{t}_{\text{ini}}
			\end{bmatrix}
			\\  0  \\ [0\;\;I_p]Y^{t}S_{\kappa,1}^{\top} +  \begin{bmatrix}
			0 & [I_p\;\; 0] y^{t}_{\text{ini}}
			\end{bmatrix}
			\end{bmatrix}, \\ 
			& U^t = \begin{bmatrix}
			U^t_p \\ U^t_f
			\end{bmatrix}, \quad Y^t = \begin{bmatrix}
			Y^t_p \\ Y^t_f
			\end{bmatrix}, \\ 
			& u_{\text{ini}}^{t+1} = S_{T_{\text{ini}},m} u_{\text{ini}}^{t} + [0\;\; I_m]^{\top} u^{\tau}_0, \\
			& y_{\text{ini}}^{t+1} = S_{T_{\text{ini}},p} y_{\text{ini}}^{t} + [0\;\; I_p]^{\top} y^{\tau}_0.
			\end{align}
		\end{subequations}
		Every $0$ above is in proper dimension. We omit the dimensions of the $0$ elements for compactness. 
	}
	Observe that the update rules \eqref{eq:sf_online} are different from~\eqref{eq:sf_solved}. The main differences between~\eqref{eq:sf_online} and~\eqref{eq:sf_solved} are that the variables $u^{\tau}$, $y^{\tau}$, $\nu^{\tau}$,  are updated through the shift matrices. \rev{In addition, $\mathcal{H}^t$ and $h^t$ are updated based on the latest input/output pair, $(u_0^{\tau},y_0^{\tau})$. }
	\rev{For $\hat{u}^{\tau} = S_{N,m}u^{\tau}$, the shift matrix $S_{N,m}$ moves $u_{k+1}^{\tau}$ in the place of $u_k^{\tau}$ for all $k = 0,\cdots,N-2$. The logic behind this updating scheme is that the prediction horizon changes, e.g., from $(t,t+1,\cdots,t+N-1)$ to $(t+1,t+2,\cdots,t+N)$, every time the optimization problem is getting updated. The updating scheme uses the shift matrices to appropriately ``initialize''  the solution of the optimization problem, $u^{\tau}$, $y^{\tau}$, $\nu^{\tau}$, at $t+1$ using the one for time $t$. The update of $\mathcal{H}^t$ and $h^t$ is to reflect the change of the optimization problem instead of the initialization of the variables. Equation~\eqref{eq:Update_H&h} simply removes the oldest input/output pair and adds the latest one, $(u_0^{\tau},y_0^{\tau})$, in the way that complies with the logistics of constructing the behavioral system model.
	}
	We summarize the ODeePC in Algorithm~\ref{alg:ODeePC}. 
	\begin{algorithm}[H]
		\caption{ODeePC 
		}
		\label{alg:ODeePC}
		\begin{algorithmic}[1]
			\State \textbf{Initialize } $\tau=t=1$, $\alpha$, $\epsilon$, $u^\tau$, $y^\tau$, $g^\tau$, $\nu^\tau$, $\Hc^{t}$, $h^{t}$
			\State \textbf{Repeat}  
			\If {mod$(\tau$,$N_I)\geq 1$}
			\State Compute~\eqref{eq:sf_solved} using Algorithm~\ref{algorithm_vector_hankel}
			\Else
			\State Apply control $u_0^{\tau}$ to the system
			\State Update $\Hc^{t+1}$ and $h^{t+1}$ \rev{by~\eqref{eq:Update_H&h} } 
			\State Compute~\eqref{eq:sf_online} using Algorithm~\ref{algorithm_vector_hankel}
			\State $t \mapsto t+1$
			\EndIf
			\State  $\tau \mapsto \tau+1$
		\end{algorithmic} 
	\end{algorithm}
	Note that Algorithm~\ref{alg:ODeePC} incorporates Algorithm~\ref{algorithm_vector_hankel}, which will be introduced in Section~\ref{sec:FFT}, to compute ${\Hc^{\tau}}^\top\nu^\tau$ and $\Hc^{\tau}g^{\tau}$ in~\eqref{eq:sf_solved} and \eqref{eq:sf_online}. The matrix-vector multiplications are the most computationally heavy part in Algorithm~\ref{alg:ODeePC}. Algorithm~\ref{algorithm_vector_hankel} gets around the multiplications with FFTs. We will first show the convergence of the ODeePC and then in Section~\ref{sec:FFT} we explain how FFT-based approach works.
	
	\subsection{Convergence of ODeePC}\label{sec:PF_convergence}
	In this section, we  show that under mild assumptions, Algorithm~\ref{alg:ODeePC} converges Q-linearly to a neighborhood of the optimal point. For convenience of notation, we denote $z^{\tau} = [{u^{\tau}}^\top \; {y^{\tau}}^\top \; {g^{\tau}}^\top \; {\nu^{\tau}}^\top]^\top$, and we rewrite~\eqref{eq:sf_online} as:
	\begin{align}\label{eq:sf_online2}
	z^{{\tau}+1} = \text{Proj}_{\mathcal{U}\times\mathcal{Y}\times \real^{\kappa+N_\nu}} \{ z^{\tau} - \alpha \Psi^{\tau}(z^{\tau}) \},
	\end{align}
	where the time-varying gradient step is embedded in $\Psi^{\tau}$.  We assume the Lipschitz continuity and monotonicity of $\Psi^{\tau}$, stated formally in Assumption~\ref{ass:Bdd_meas_err}. 
	\begin{assumption}\longthmtitle{Lipschitz continuity and monotonicity of the gradient}\label{ass:Bdd_meas_err}
		There exists a finite constant $\sigma_{\Psi}\in\real_+$ such that $\norm{\Psi^{\tau}(z^1) - \Psi^{\tau}(z^2)}\leq \sigma_{\Psi}\norm{z^1 - z^2}$ for all $z^1, z^2 \in \mathcal{U}\times\mathcal{Y}\times \real^{\kappa+N_\nu}$ for all $\tau$. In addition,  $\Psi^{\tau}$ is strongly monotone with constant $\eta$.
	\end{assumption}
	Assumption~\ref{ass:Bdd_meas_err} is actually equivalent to \cite[Lemma 3.4]{koshal2011multiuser}.  
	The reason of making Assumption~\ref{ass:Bdd_meas_err} instead of deriving~\cite[Lemma 3.4]{koshal2011multiuser} for our case is that the proof of~\cite[Lemma 3.4]{koshal2011multiuser} will require breaking the equality constraint on $g$ into inequality constraints. In addition, because the gradient mapping $\Psi^{\tau}$ has two different forms (dependent on the iteration step), respectively \eqref{eq:sf_solved} and \eqref{eq:sf_online}, we need to show the Lipschitz continuity and strong monotonicity for both forms. 
	Those are tedious and not the focus of this paper. 
	
	We further make Assumption~\ref{ass:Bdd_opt} that considers that the variation of the optimal point is bounded, enabling each time close tracking of the optimal point within a reasonable number of iterations of ODeePC.
	\begin{assumption}\longthmtitle{Bounded variation of the optimal point}\label{ass:Bdd_opt}
		There exists a finite constant $\sigma_z\in\real_+$ such that the optimal points for the consecutive time steps $t$ (or $\tau$) and $t+1$ (or $\tau+N_I$) of optimization~\eqref{eq:ODeePC} satisfy $\norm{z^{\tau+N_I,\star} - z^{\tau,\star}}_2\leq \sigma_z$. 
	\end{assumption}
	Recall that we index $z$ by the iteration number rather than the time instances $t$. The optimal point of $z$ for iterations $\tau,\cdots,\tau+N_I-1$ (or time $t$) stays unchanged because the associated optimization remains the same. The optimal point changes for iteration $\tau+N_I$ (time $t+1$). 
	\rev{
		We last make Assumption~\ref{ass:persist_exciting} before showing the convergence results of ODeePC.
		\begin{assumption}\longthmtitle{Persistently exciting condition for online updates} \label{ass:persist_exciting}
			For every time $t\in\mathbb{N}$, $U^t$ has full row rank.
		\end{assumption}
		Assumption~\ref{ass:persist_exciting} ensures that for every time instance, the sequence of the controls used for constructing the block Hankel matrix, $U^t$, satisfies the persistently exciting condition. If the persistently exciting condition is lost throughout the course of updates of the block Hankel matrices, then the optima of~\eqref{eq:ODeePC} may not correspond to optimal control for the underlying control system. It is very challenging to ensure the persistently exciting condition being satisfied for all times (especially in an online setup). Specifically, if the control input stays close to certain value for a period of time, it can destroy the richness of the data and, as a result, the persistently exciting condition may no longer hold. To lower the odds of the violation of this condition, one may halt updating the block Hankel matrices once the control input (approximately) reaches a  steady state. A comprehensive method to avoid the violation of this condition is a challenging future work. We make it as an assumption in this paper.
	}
	
	With all the elements in place, we state Theorem~\ref{convergence} which establishes convergence of Algorithm~\ref{alg:ODeePC}.
	\begin{theorem}\longthmtitle{Convergence of ODeePC}
		\label{convergence}
		\footnote{We assume that the measurement noise is zero in the proof. }
		If Assumptions~\ref{ass:Bdd_meas_err}, \ref{ass:Bdd_opt} and~\rev{\ref{ass:persist_exciting}} hold, and $\rho(\alpha):= \sqrt{1+\alpha^2\sigma_{\Psi}^2-2\alpha\eta}<1$, then Algorithm~\ref{alg:ODeePC} has $z^\tau$ converge Q-linearly to a  neighborhood of optimal point of~\eqref{eq:ODeePC}, given as:
		\begin{align*}
		\limsup_{\tau\rightarrow \infty} \|z^{\tau} - z^{\tau,\star}\|_2 =  \frac{\sigma_z}{1-\rho(\alpha)}.
		\end{align*}
	\end{theorem}
	The proof of Theorem~\ref{convergence} is similar to the proof of \cite[Theorem 4]{bernstein2019real}. The main element of the proof is the derivation of an inequality between $\|z^{\tau+1} - z^{\tau+1,\star}\|$ and $\|z^{\tau} - z^{\tau,\star}\|$. In the proof, Assumption~\ref{ass:Bdd_meas_err} is used for substituting the terms $(\Psi^{\tau}(z^{\tau}) - \Psi^{\tau}(z^{\tau,\star}))^\top (z^{\tau} - z^{\tau,\star})$ and $\|\Psi^{\tau}(z^{\tau}) - \Psi^{\tau}(z^{\tau,\star})\|_2^2$ by the product of some coefficient (determined by $\sigma_{\Psi}$ and $\eta$) and $\|z^{\tau} - z^{\tau,\star}\|$. Readers are referred to \cite[Theorem 4]{bernstein2019real} for more details. 

	\subsection{Fast Computation of Hankel Matrix-Vector Product}\label{sec:FFT}
	In this section, we describe Algorithm~\ref{algorithm_vector_hankel}, which efficiently computes the matrix-vector multiplications in~\eqref{eq:sf_solved} and~\eqref{eq:sf_online}. We will use the fact that $\Hc$ in~\eqref{eq:sf_solved} and~\eqref{eq:sf_online} is the concatenation of block Hankel matrices $U$ and $Y$, and we apply FFT to exploit the convolutional structure of the Hankel matrices.
	
	We slightly abuse the notation in this subsection, and we let $n$ and $m$ denote the dimensions of a general Hankel matrix $\mathbf{H}$ so that $\mathbf{H}\in\mathbb{R}^{n\times m}$, unrelated to the state-space model representation \eqref{eq:dss}. Given that, the main problem we are concerned with here can be described as follows. 
	\begin{problem}
		Given a $m$-element vector $v\in\mathbb{R}^{m}$:
		\begin{align}
		v:=\begin{pmatrix} v_1 & v_2 & \cdots & v_m \end{pmatrix}^\top, \label{v}
		\end{align}
		and a Hankel matrix $\mathbf{H}\in\mathbb{R}^{n\times m}$, where $h_i\in\mathbb{R},\;\forall i$, defined as:
		\begin{align}
		\mathbf{H}:=\begin{pmatrix}
		h_1 & h_2 & \cdots & h_{m-1} & h_m\\
		h_2 & h_3 & \cdots & h_{m} & h_{m+1}\\
		\vdots & \vdots \\
		h_{n-1} & h_n & \cdots & h_{n+m-3} & h_{n+m-2}\\
		h_n & h_{n+1} & \cdots & h_{n+m-2} & h_{n+m-1}\\
		\end{pmatrix},\hspace{3mm} \label{Hankel}
		\end{align}
		compute the product $p_{\text{H}}=\mathbf{H}v$ in a computationally efficient manner.
	\end{problem}
	Building on the algorithm proposed in \cite{FFTalgorithm} for computing the product of a square Hankel matrix and a vector, we propose Algorithm~\ref{algorithm1} for computing the product $p_{\text{H}}=\mathbf{H}v$ efficiently for \textit{non-square Hankel matrices}. The complexity of Algorithm~\ref{algorithm1} is introduced through the following theorem.
	
	\begin{algorithm}
		\caption{Fast Hankel matrix-vector product}
		Given a vector $v\in\mathbb{R}^m$ and a Hankel matrix $\mathbf{H}\in\mathbb{R}^{n\times m}$, compute the vector $p_{\text{H}}=\mathbf{H}v$ via the following steps.
		\begin{enumerate}
			\item Define a new $(n+m-1)$-element vector $c$ as:
			\begin{align*}
			c=\begin{pmatrix}h_m &  h_{m+1}  &  \cdots &  h_{n+m-1} &  h_1 &  h_2 &  \cdots &  h_{m-1} \end{pmatrix}^\top
			\end{align*}	
			\item Define a $(n+m-1)$-element vector $v_e\in\mathbb{R}^{n+m-1}$ by permuting the vector $v$ and adding $(n-1)$ zeros so that:
			\begin{align*}
			v_e=\begin{pmatrix} v_m & v_{m-1} & \cdots & v_1 &  0 & \cdots & 0 \end{pmatrix}^\top
			\end{align*}
			\item Compute a $(n+m-1)$-element vector $y$ as:
			\begin{align*}
			y=\mathbf{IFFT}(\mathbf{FFT}(c)\circ\mathbf{FFT}(v_e))
			\end{align*}
			where $(\circ)$ is the Hadamard product of the two vectors, $\mathbf{FFT}$ is the fast Fourier transform, and $\mathbf{IFFT}$ is the inverse fast Fourier transform.
			\item Let $y=\begin{pmatrix}y_1 & y_2 & \cdots & y_{(n+m-2)} & y_{(n+m-1)}\end{pmatrix}^\top, y\in\mathbb{R}^{n+m-1}$. Then the product $p_{\text{H}}=\mathbf{H}v$ is given by:
			\begin{align*}
			p_{\text{H}}=\begin{pmatrix} y_1 & y_2 & \cdots & y_{n-1} & y_n \end{pmatrix}^\top
			\end{align*}
			i.e., the subvector defined by the first $n$ elements of the vector $y$.  
		\end{enumerate}
		\label{algorithm1}	
	\end{algorithm}
	
	\begin{theorem}
		The complexity of Algorithm~\ref{algorithm1} that computes the product of a Hankel matrix $\mathbf{H}\in\mathbb{R}^{n\times m}$ and a vector $v\in\mathbb{R}^m$ is $O\big(\max(n,m)\log(\max(n,m)) \big)$.
	\end{theorem}
	\begin{pf}
		Following the same steps as in \cite{FFTalgorithm}, we first note that each FFT has a complexity $5(n+m-1)\log(n+m-1)$, and the pointwise multiplication has a complexity $6(n+m-1)$ and the inverse FFT $5(n+m-1)\log(n+m-1)$. By combining all of these together, we have that the complexity of the algorithm is:
		\begin{align}\label{eq:N_multiplication}
		15(n+m-1)\log(n+m-1)+6(n+m-1).
		\end{align}
		Equation~\eqref{eq:N_multiplication} can be rewritten in Big-O notation, which is the desired $O\big(\max(n,m)\log(\max(n,m)) \big)$.
	\end{pf}
	It remains to show that the proposed algorithm precisely computes the desired product $p_{\text{H}}=\mathbf{H}v$. This is carried out through the following theorem.
	\begin{theorem}
		Consider a Hankel matrix $\mathbf{H}\in\mathbb{R}^{n\times m}$ as defined in \eqref{Hankel} and a vector $v\in\mathbb{R}^m$ as defined in \eqref{v}. The first $n$ elements of the vector $y$ obtained from the following computation:
		\begin{align}
		y=\mathbf{IFFT}(\mathbf{FFT}(c)\circ\mathbf{FFT}(v_e))
		\end{align}
		where $c$ and $v_e$ are $(n+m-1)$-element vectors, defined as:
		\begin{align*}
		c&:=\begin{pmatrix}h_m &   h_{m+1}  &  \cdots &  h_{n+m-1} &  h_1 &  h_2 &  \cdots &  h_{m-1} \end{pmatrix}^\top, 
		\\
		v_e&:=\begin{pmatrix} v_m & v_{m-1} & \cdots & v_1 &  0 & \cdots & 0 \end{pmatrix}^\top, 
		\end{align*}
		yields the exact product $p_{\text{H}}=\mathbf{H} v$.
		\label{Hankelvectorprod}
	\end{theorem}
	\begin{pf}
		The proof is deferred to the appendix.
	\end{pf}
	Through this theorem, we have established that Algorithm~\ref{algorithm1} computes the product $p_{\text{H}}=\mathbf{H}v$. By building on these results, we will now design an algorithm that will allows us to efficiently compute the product of a \textit{block} Hankel matrix $\mathbf{H}$ (comprising column vectors $\mathbf{h}_i$) and a vector $v$.
	The main problem we seek to address can be stated as follows. 
	\begin{problem}
		Given a $m$-element vector $v\in\mathbb{R}^m$:
		\begin{align}
		v:=\begin{pmatrix} v_1 & v_2 & \cdots & v_m \end{pmatrix}^T \label{v_vec},
		\end{align}
		and a block Hankel matrix $\mathbf{H}\in\mathbb{R}^{(nl)\times m}$, where $\mathbf{h}_i\in\mathbb{R}^l,\;l>1,\;\forall i$, defined as:
		\begin{align}
		\mathbf{H}:=\begin{pmatrix}
		\mathbf{h}_1 & \mathbf{h}_2 & \cdots & \mathbf{h}_{m-1} & \mathbf{h}_m\\
		\mathbf{h}_2 & \mathbf{h}_3 & \cdots & \mathbf{h}_{m} & \mathbf{h}_{m+1}\\
		\vdots & \vdots \\
		\mathbf{h}_{n-1} & \mathbf{h}_n & \cdots & \mathbf{h}_{n+m-3} & \mathbf{h}_{n+m-2}\\
		\mathbf{h}_n & \mathbf{h}_{n+1} & \cdots & \mathbf{h}_{n+m-2} & \mathbf{h}_{n+m-1}\\
		\end{pmatrix},\hspace{3mm} \label{Hankel_vec}
		\end{align}
		Compute the product $p_{\text{BH}}=\mathbf{H}v$ in a computationally efficient manner.
		\label{Problem_vec_Hankel}
	\end{problem}
	Our goal is to compute the product of a block Hankel matrix $\mathbf{H}$ that comprises column vectors $\mathbf{h}_i\in\mathbb{R}^l$, with a vector $v$. We propose the following algorithm for solving Problem~\eqref{Problem_vec_Hankel}.
	
	\begin{algorithm}
		\caption{Fast block Hankel matrix - vector product}
		The product $p_{\text{BH}}$ of a vector $v\in\mathbb{R}^m$ and a block Hankel matrix $\mathbf{H}\in\mathbb{R}^{(nl)\times m}$ can be computed through the following steps.
		\begin{enumerate}
			\item \textbf{Initialize} $\beta=1$ and define the vector: $c\in\mathbb{R}^{l(n+m-1)}$
			\begin{align*}
			c=\begin{pmatrix}\mathbf{h}_m^\top   & \cdots & \mathbf{h}_{n+m-1}^\top & \mathbf{h}_1^\top & \mathbf{h}_2^\top & \cdots & \mathbf{h}_{m-1}^\top \end{pmatrix}^\top
			\end{align*}	
			\item Define a $(n+m-1)$-element vector $v_e\in\mathbb{R}^{n+m-1}$ by permuting the vector $v$ and adding $(n-1)$ zeros so that:
			\begin{align*}
			v_e=\begin{pmatrix} v_m & v_{m-1} & \cdots & v_1 &  0 & \cdots & 0 \end{pmatrix}^T 
			\end{align*}
			\item \label{alg:step3} While $\beta\leq l$, construct:
			\begin{align*}
			c^{(\beta)}=\begin{pmatrix}\mathbf{h}_m(\beta)  &  \cdots &  \mathbf{h}_{n+m-1}(\beta) &  \mathbf{h}_1(\beta)  &  \cdots &  \mathbf{h}_{m-1}(\beta) \end{pmatrix}^\top
			\end{align*}	
			and apply \textbf{Algorithm}~\ref{algorithm1}, where $c=c^{(\beta)}$, to obtain $y^{(\beta)}=y\in\mathbb{R}^{n+m-1}$. Update $\beta \rightarrow \beta+1$.
			\item Compute the product $p_{\text{BH}}=\mathbf{H}v\in\mathbf{R}^{n l}$ as:
			\begin{align*}
			p_{\text{BH}}=\begin{pmatrix} y^{(1)}(1)  &  \cdots &  y^{(l)}(1)  &  \cdots &  y^{(1)}(n)  &  \cdots &  y^{(l)}(n) \end{pmatrix}^\top 
			\end{align*}
		\end{enumerate}
		\label{algorithm_vector_hankel}	
	\end{algorithm}
	Following the steps in the proof of Theorem~\ref{Hankelvectorprod}, one can show that application of Algorithm~\ref{algorithm_vector_hankel} indeed results in the precise computation of the product $p_{\text{BH}}$ of a block Hankel matrix $\mathbf{H}$ and a vector $v$. 
	
	
	\section{Numerical Validation of ODeePC}
	\label{sec:ODeePCVal}
	We consider the following predictive control problem with a linear time-varying dynamic system:
	\begin{align}\label{eq:sim_model}
	& \underset{x, u\in\mathcal{U}, \rev{y}}{\text{minimize}}
	\sum_{k=0}^{N-1}\norm{y_k-r_{t+k}}^2_2,\\ \nonumber
	& \text{s. t. }
	\hspace{5mm} x_{k+1} = A_{t+k}x_k + B_{t+k}u_k, \forall k \in \{0, \cdots, N-1 \}, \\ \nonumber
	& \hspace{12mm} y_k = Cx_k, \forall k \in \{0, \cdots, N-1 \},  \\ \nonumber
	&\hspace{12mm}  x_0=\hat{x}_t,
	\end{align}
	where \rev{$\mathcal{U} = \{u \vert \norm{u}_{\infty}\leq 1\}$}, $A_{t}\in\real^{10\times 10}$ and  $B_{t}\in\real^{10\times 10}$ for every time instance $t$, $C\in\real^{10\times 10}$, and $r_t$ is the reference signal. 
	\rev{The objective of the predictive control is to track a given reference output signal $r_t$. All $A_{t}$, $B_{t}$ and $C$ are randomly generated with $\norm{A_{t}}_2 = \norm{B_{t}}_2 = \norm{C}_2 = 1$. The generated system was verified controllable and observable. Both $A_t$ and $B_t$ vary with $t$, specifically, the magnitude of each entry of $A_t$ increases $p\%$ of its value whenever $t$ changes to $t+1$, where $p$ is randomly generated with a uniformly distribution in $[-0.01,0.01]$. Similar conditions are imposed on $B_t$.
	}%
	The references $r_t$ are changed randomly with uniform distribution between $[0, 0.1]$ for every $1000$ iterations of $t$. 
	
	To apply DeePC or ODeePC to the predictive control problem defined by~\eqref{eq:sim_model}, we construct the data-based model~\eqref{eq:DeePC_model} by prerunning a number of iterations with a sequence of $u_t$. Optimization~\eqref{eq:sim_model} can then be reformulated in the form of~\eqref{eq:ODeePC} (the objective function can be time-varying without loss of generality). \rev{Analytical results for LTI systems~\cite[Lemma 1]{DeePC} suggest that $T_{\text{ini}}$ larger than the dimension of $x$, which is $10$, is sufficient to observe the system behavior. We use this as a reference and choose $T_{\text{ini}} = 20$ for this LTV system. The initial data also satisfy persistently exciting condition of order $10(T_{\text{ini}}+N)$.} The parameters of the optimization and model are given in Table~\ref{table:numbers}.
	\begin{table}[H]
		\centering 
		\caption{Number setup for the simulations.}
		\begin{tabular}{|c|c|c|c|c|} \hline
			$N_I$ & $T_{ini}$ & $N$ &$\kappa$ & $\epsilon_g$ \\ \hline 
			50 & 20 & \rev{120} & 1651 & 0.1 \\ \hline
		\end{tabular}
		\label{table:numbers} 
	\end{table}
	
\begin{figure}[H]
		\centering
		\includegraphics[width=1\columnwidth]{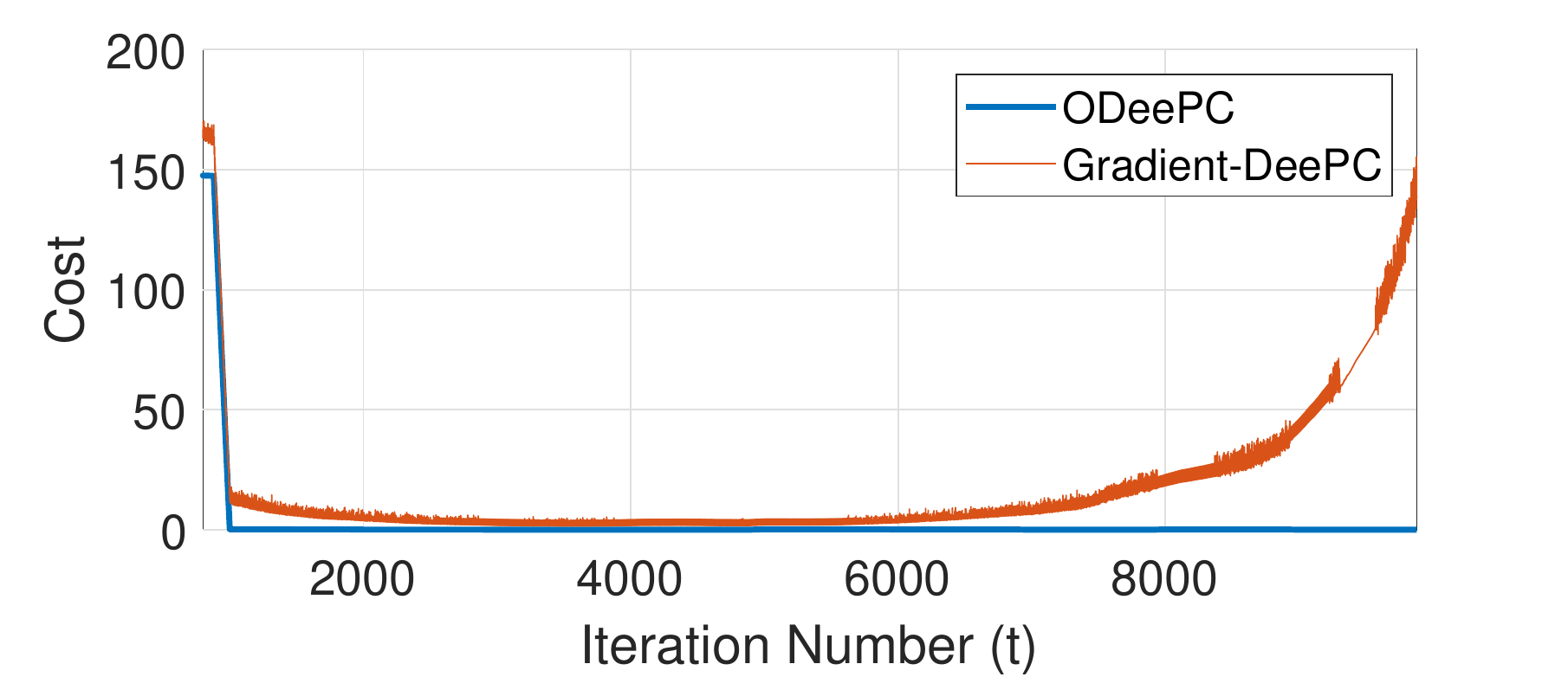}
		\caption{Evolution of costs over iterations for ODeePC and Gradient-DeePC algorithms.}
		\label{fig:Cost}
	\end{figure}

\begin{figure}[H]
		\centering
		\includegraphics[width=1\columnwidth]{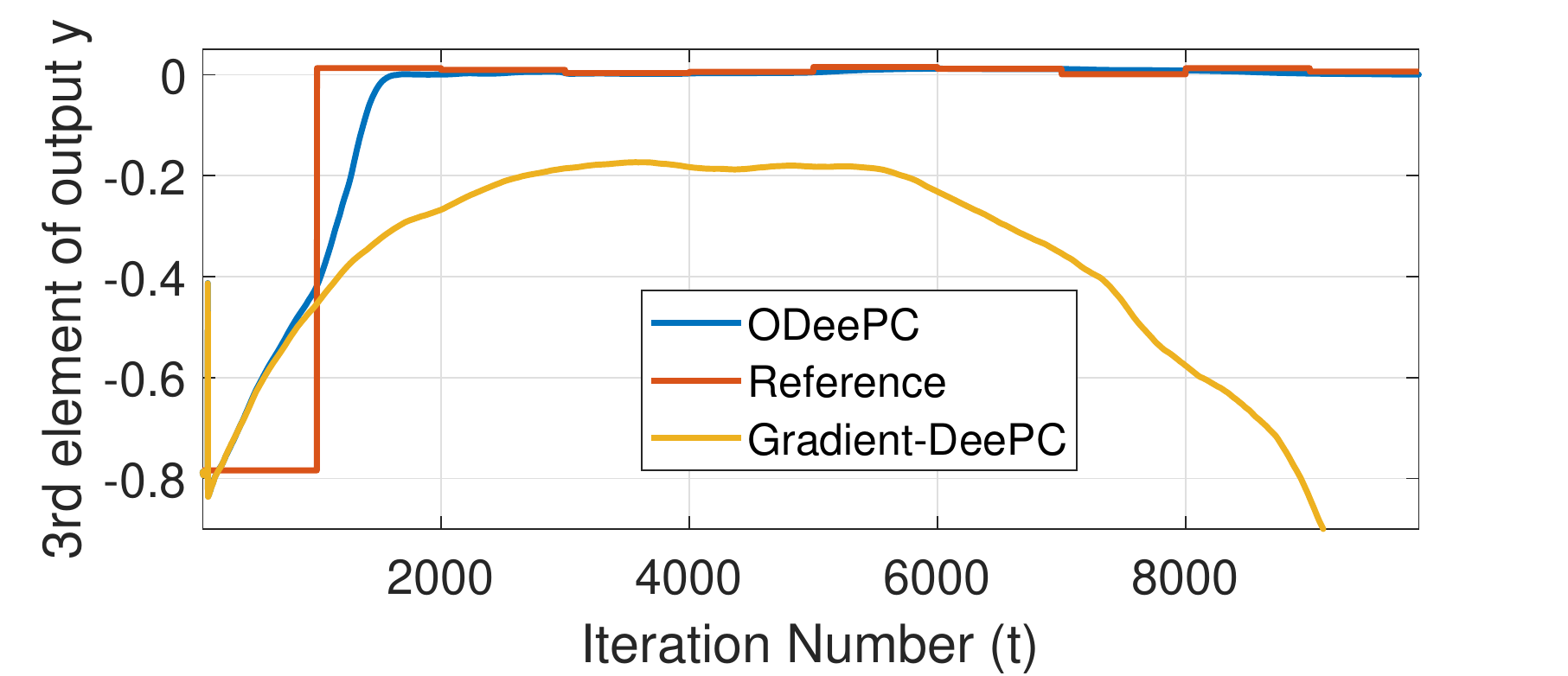}
		\caption{Tracking performance of one entry of the output $y$.}
		\label{fig:Tracking-03}
	\end{figure}

\begin{figure}[H]
		\centering
		\includegraphics[width=1\columnwidth]{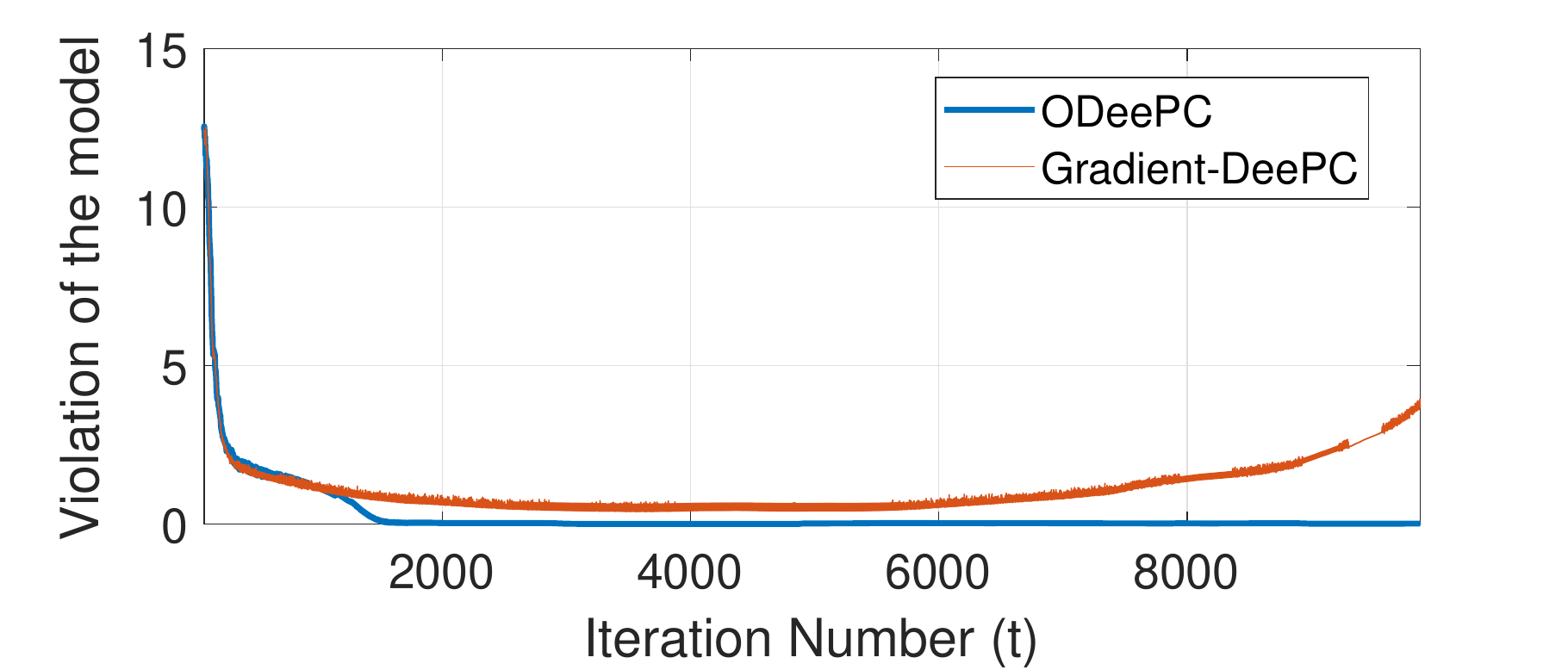}
		\caption{Two-norm violation of equality~\eqref{eq:DeePC_model}.}
		\label{fig:Hg-h-A00001.eps}
	\end{figure}

	
	\rev{In this numerical study, we first apply  DeePC to solve~\eqref{eq:sim_model}. However, because of the time-varying $A_t$ and $B_t$, the optimization associated with DeePC becomes infeasible in less than 10 iterations. We therefore switch to an alternative gradient-DeePC algorithm which is the same as Algorithm~\ref{alg:ODeePC} except that $\Hc^t$ is kept unchanged over time.} We view the gradient-DeePC as a varient of DeePC (or benchmark) to compare the performance of ODeePC. The results are shown in Figures \ref{fig:Cost}, \ref{fig:Tracking-03}, and~\ref{fig:Hg-h-A00001.eps}. It can be seen that Gradient-DeePC starts diverging at around $t= 6000$; on the contrary, ODeePC  converges to a near optimal point with the cost close to zero.  
	
	
	Another reason for the performance difference between ODeePC and Gradient-DeePC is on the element shifting step embedded in~\eqref{eq:sf_online} when the optimization is updated from $t$ to $t+1$. Recall that ODeePC updates both $\Hc$ and $h$, whereas Gradient-DeePC updates $h$ and keeps $\Hc$ unchanged. This makes a difference because the shifting on variable $z$  results in proper initialization of the optimization problem at $t+1$ using the solution at $t$. For the case of ODeePC, the majority of the constraints defined by $\Hc g = h$ are kept when $t$ is changed to $t+1$, so the shifts defined in~\eqref{eq:sf_online} result in variables that satisfy most equality constraints of the optimization problem at $t+1$. On the contrary, Gradient-DeePC only shifts $h$ when the optimization problem is updated. This implies that all the constraints of the optimization problem at $t+1$ are different from the ones at $t$. The ``disconnection'' between consecutive instances of the optimization problem is an additional disadvantage of Gradient-DeePC compared to ODeePC, which partially causes the divergence observed in Fig.~\ref{fig:Hg-h-A00001.eps}.
	
	
	We also compare the computational times between Algorithm~\ref{algorithm_vector_hankel} and the direct matrix vector multiplication method. In our rather moderate-scale example, Algorithm~\ref{algorithm_vector_hankel} takes $0.51$ second on average with $N_I= 50$ inner gradient iterations, whereas the average time for the direct multiplication is approximately $1.1$ seconds. Both simulations are done on a desktop with 3.5-GHz CPU and 16-GB RAM. Algorithm~\ref{algorithm_vector_hankel} slashes the computational time in half. 
	
	\rev{
		To further explore the benefits of  Algorithm~\ref{algorithm_vector_hankel} in real-time ODeePC implementation, we simulate a larger system with $A_t\in\real^{80\times 80}$, $B_t\in\real^{80\times 80}$, and $C\in\real^{60\times 80}$. The variables $T_{\text{ini}}$ and $N$ are adjusted accordingly with $T_{\text{ini}} = N = 80$. We note that $N$ can be chosen freely for computational time comparison purpose because the prediction horizon solely depends on the application. A noticeable improvement is found: Algorithm~\ref{algorithm_vector_hankel} took  $219.94$ seconds compared to $1040.73$ seconds for direct multiplication method ($N_I=50$). Compared to the original test case with $m = p = 10$, the computational time grows significantly in this case with $m = 80$ and $p = 60$ (from $1.1$ to $1043.73$ seconds). The reason behind it is that the size of the block Hankel matrices roughly grows quadratically with respect to $m$ and the dimension of $x$, $n$. To see this, note that sufficiently large $T_{\text{ini}}$ is at the same scale of $n$. In addition, the persistently exciting condition requires row rank of the block Hankel matrix being $m \cdot (T_{\text{ini}} +N) \approx m (n +N)$. To satisfy the persistently exciting condition, we need data to construct block Hankel matrix with number of columns not less than $m \cdot (n+N)$. This translates to the size of the block Hankel matrix approximately at $\big(m(n +N)\big) \times \big(m(n +N)\big)$, and, as a result, the number of multiplications in~\eqref{eq:sf_online} is approximated at the order of $\big(m(n +N) \big)^2$. Thus, online implementation of ODeePC needs to be cautious regarding the complexity growth originated from the data.
		
	}
	
	
	\section{Conclusion and Future Work}
	\label{sec:conc}
	In this paper, we presented an online data-enabled predictive (ODeePC) control method for optimal control of unknown systems, building on the recently proposed DeePC \cite{DeePC}. Our proposed ODeePC method leverages a primal-dual algorithm with real-time measurement feedback and recorded system data to compute the optimal control policy in real-time as system conditions change. ODeePC can generate control inputs dynamically, tracking changes in the system operating point while manifesting high computational efficiency. We prove that ODeePC's iterative update rule converges to a neighborhood of the optimal control policy. 
	
	\rev{Though ODeePC provides an interesting model-free option for online predictive control, we would like to note that standard methods that combine online parameter estimation and MPC might outperform ODeePC if system model structure (linear or certain class of nonlinear function) is known a-priori. One reason is that the standard method implicitly incorporates the model structure while ODeePC does not. Another reason is that the standard method records system behavior by latest parameters with dimension much lower than the collection of data in ODeePC. However, when it comes to control of black box system so that parameter estimation approaches become less effective, ODeePC provides a good heuristic. In this context, DeePC is similar conceptually to system identification plus MPC; whereas ODeePC is similar to online parameter identification plus MPC. In-depth comparison between those comparable concepts with different approaches is considered as an important future work. We will also attempt to find a control strategy that preserves persistently exciting condition for online update of the behavioral system model. 
		
	}
	
	\bibliographystyle{IEEEtran}
	\bibliography{bib_DeePC}

\begin{thebibliography}{10}
\providecommand{\url}[1]{#1}
\csname url@samestyle\endcsname
\providecommand{\newblock}{\relax}
\providecommand{\bibinfo}[2]{#2}
\providecommand{\BIBentrySTDinterwordspacing}{\spaceskip=0pt\relax}
\providecommand{\BIBentryALTinterwordstretchfactor}{4}
\providecommand{\BIBentryALTinterwordspacing}{\spaceskip=\fontdimen2\font plus
\BIBentryALTinterwordstretchfactor\fontdimen3\font minus
  \fontdimen4\font\relax}
\providecommand{\BIBforeignlanguage}[2]{{%
\expandafter\ifx\csname l@#1\endcsname\relax
\typeout{** WARNING: IEEEtran.bst: No hyphenation pattern has been}%
\typeout{** loaded for the language `#1'. Using the pattern for}%
\typeout{** the default language instead.}%
\else
\language=\csname l@#1\endcsname
\fi
#2}}
\providecommand{\BIBdecl}{\relax}
\BIBdecl

\bibitem{DeePC}
J.~Coulson, J.~Lygeros, and F.~D{\"o}rfler, ``Data-enabled predictive control:
  In the shallows of the {D}ee{PC},'' in \emph{2019 18th European Control
  Conference (ECC)}.\hskip 1em plus 0.5em minus 0.4em\relax IEEE, 2019, pp.
  307--312.

\bibitem{Baros}
A.~M. Annaswamy, A.~R. Malekpour, and S.~Baros, ``Emerging research topics in
  control for smart infrastructures,'' \emph{Annual Reviews in Control},
  vol.~42, pp. 259--270, 2016.

\bibitem{Dorfler_4}
J.~B. Rawlings and D.~Q. Mayne, \emph{Model predictive control: Theory and
  design}.\hskip 1em plus 0.5em minus 0.4em\relax Nob Hill Publishing, 2009.

\bibitem{Dorfler_5}
E.~F. Camacho and C.~B. Alba, \emph{Model predictive control}.\hskip 1em plus
  0.5em minus 0.4em\relax Springer Science \& Business Media, 2013.

\bibitem{Dorfler_6}
A.~Bemporad and M.~Morari, ``Robust model predictive control: A survey,'' in
  \emph{Robustness in identification and control}.\hskip 1em plus 0.5em minus
  0.4em\relax Springer, 1999, pp. 207--226.

\bibitem{Dorfler_7}
D.~Q. Mayne, ``Model predictive control: Recent developments and future
  promise,'' \emph{Automatica}, vol.~50, no.~12, pp. 2967--2986, 2014.

\bibitem{Dorfler_8}
F.~Borrelli, A.~Bemporad, and M.~Morari, \emph{Predictive control for linear
  and hybrid systems}.\hskip 1em plus 0.5em minus 0.4em\relax Cambridge
  University Press, 2017.

\bibitem{Dorfler_1}
F.~Lamnabhi-Lagarrigue, A.~Annaswamy, S.~Engell, A.~Isaksson, P.~Khargonekar,
  R.~M. Murray, H.~Nijmeijer, T.~Samad, D.~Tilbury, and P.~Van~den Hof,
  ``Systems \& control for the future of humanity, research agenda: Current and
  future roles, impact and grand challenges,'' \emph{Annual Reviews in
  Control}, vol.~43, pp. 1--64, 2017.

\bibitem{Dorfler_2}
Z.-S. Hou and Z.~Wang, ``From model-based control to data-driven control:
  Survey, classification and perspective,'' \emph{Information Sciences}, vol.
  235, pp. 3--35, 2013.

\bibitem{Dorfler_22}
M.~C. Campi and E.~Weyer, ``Finite sample properties of system identification
  methods,'' \emph{IEEE Transactions on Automatic Control}, vol.~47, no.~8, pp.
  1329--1334, 2002.

\bibitem{Dorfler_23}
M.~Vidyasagar and R.~L. Karandikar, ``A learning theory approach to system
  identification and stochastic adaptive control,'' in \emph{Probabilistic and
  randomized methods for design under uncertainty}.\hskip 1em plus 0.5em minus
  0.4em\relax Springer, 2006, pp. 265--302.

\bibitem{Dorfler_24}
S.~Tu, R.~Boczar, A.~Packard, and B.~Recht, ``Non-asymptotic analysis of robust
  control from coarse-grained identification,'' \emph{arXiv preprint
  arXiv:1707.04791}, 2017.

\bibitem{Dorfler_16}
F.~L. Lewis, D.~Vrabie, and K.~G. Vamvoudakis, ``Reinforcement learning and
  feedback control: Using natural decision methods to design optimal adaptive
  controllers,'' \emph{IEEE Control Systems Magazine}, vol.~32, no.~6, pp.
  76--105, 2012.

\bibitem{Dorfler_17}
Y.~Ouyang, M.~Gagrani, and R.~Jain, ``Learning-based control of unknown linear
  systems with thompson sampling,'' \emph{arXiv preprint arXiv:1709.04047},
  2017.

\bibitem{Dorfler_18}
B.~Kiumarsi, F.~L. Lewis, H.~Modares, A.~Karimpour, and M.-B. Naghibi-Sistani,
  ``Reinforcement {Q}-learning for optimal tracking control of linear
  discrete-time systems with unknown dynamics,'' \emph{Automatica}, vol.~50,
  no.~4, pp. 1167--1175, 2014.

\bibitem{Dorfler_19}
A.~M. Devraj and S.~Meyn, ``Zap q-learning,'' in \emph{Advances in Neural
  Information Processing Systems}, 2017, pp. 2235--2244.

\bibitem{Dorfler_20}
B.~Recht, ``A tour of reinforcement learning: The view from continuous
  control,'' \emph{Annual Review of Control, Robotics, and Autonomous Systems},
  vol.~2, pp. 253--279, 2019.

\bibitem{Dorfler_21}
R.~Islam, P.~Henderson, M.~Gomrokchi, and D.~Precup, ``Reproducibility of
  benchmarked deep reinforcement learning tasks for continuous control,''
  \emph{arXiv preprint arXiv:1708.04133}, 2017.

\bibitem{Dorfler_28}
F.~Berkenkamp, M.~Turchetta, A.~Schoellig, and A.~Krause, ``Safe model-based
  reinforcement learning with stability guarantees,'' in \emph{Advances in
  neural information processing systems}, 2017, pp. 908--918.

\bibitem{Dorfler_29}
J.~F. Fisac, A.~K. Akametalu, M.~N. Zeilinger, S.~Kaynama, J.~Gillula, and
  C.~J. Tomlin, ``A general safety framework for learning-based control in
  uncertain robotic systems,'' \emph{IEEE Transactions on Automatic Control},
  vol.~64, no.~7, pp. 2737--2752, 2018.

\bibitem{willems1986time_1}
J.~C. Willems, ``From time series to linear system - part {I}. finite
  dimensional linear time invariant systems,'' \emph{Automatica}, vol.~22,
  no.~5, pp. 561--580, 1986.

\bibitem{willems1986time_2}
------, ``From time series to linear system - part {II}. exact modelling,''
  \emph{Automatica}, vol.~22, no.~6, pp. 675--694, 1986.

\bibitem{willems1987time_3}
------, ``From time series to linear system—part {III}: Approximate
  modelling,'' \emph{Automatica}, vol.~23, no.~1, pp. 87--115, 1987.

\bibitem{willems2005note}
J.~C. Willems, P.~Rapisarda, I.~Markovsky, and B.~L. De~Moor, ``A note on
  persistency of excitation,'' \emph{Systems \& Control Letters}, vol.~54,
  no.~4, pp. 325--329, 2005.

\bibitem{markovsky2006exact}
I.~Markovsky, J.~C. Willems, S.~Van~Huffel, and B.~De~Moor, \emph{Exact and
  approximate modeling of linear systems: A behavioral approach}.\hskip 1em
  plus 0.5em minus 0.4em\relax SIAM, 2006.

\bibitem{favoreel1999model}
W.~Favoreel, B.~De~Moor, P.~Van~Overschee, and M.~Gevers, ``Model-free
  subspace-based lqg-design,'' in \emph{Proceedings of the 1999 American
  Control Conference}, vol.~5, 1999, pp. 3372--3376.

\bibitem{kadali2003data}
R.~Kadali, B.~Huang, and A.~Rossiter, ``A data driven subspace approach to
  predictive controller design,'' \emph{Control engineering practice}, vol.~11,
  no.~3, pp. 261--278, 2003.

\bibitem{markovsky2008data}
I.~Markovsky and P.~Rapisarda, ``Data-driven simulation and control,''
  \emph{International Journal of Control}, vol.~81, no.~12, pp. 1946--1959,
  2008.

\bibitem{AndreyOnlineOpt}
A.~Bernstein, E.~Dall'Anese, and A.~Simonetto, ``Online primal-dual methods
  with measurement feedback for time-varying convex optimization,'' \emph{IEEE
  Transactions on Signal Processing}, vol.~67, no.~8, pp. 1978--1991, 2019.

\bibitem{bernstein2019real}
A.~Bernstein and E.~Dall’Anese, ``Real-time feedback-based optimization of
  distribution grids: A unified approach,'' \emph{IEEE Transactions on Control
  of Network Systems}, vol.~6, no.~3, pp. 1197--1209, 2019.

\bibitem{behaviorsystemstheory}
I.~Markovsky and P.~Rapisarda, ``Data-driven simulation and control,''
  \emph{International Journal of Control}, vol.~81, no.~12, pp. 1946--1959,
  2008.

\bibitem{Dorfler_9}
M.~Brown, J.~Funke, S.~Erlien, and J.~C. Gerdes, ``Safe driving envelopes for
  path tracking in autonomous vehicles,'' \emph{Control Engineering Practice},
  vol.~61, pp. 307--316, 2017.

\bibitem{Dorfler_10}
I.~Prodan, S.~Olaru, R.~Bencatel, J.~B. de~Sousa, C.~Stoica, and S.-I.
  Niculescu, ``Receding horizon flight control for trajectory tracking of
  autonomous aerial vehicles,'' \emph{Control Engineering Practice}, vol.~21,
  no.~10, pp. 1334--1349, 2013.

\bibitem{koshal2011multiuser}
J.~Koshal, A.~Nedi{\'c}, and U.~V. Shanbhag, ``Multiuser optimization:
  Distributed algorithms and error analysis,'' \emph{SIAM Journal on
  Optimization}, vol.~21, no.~3, pp. 1046--1081, 2011.

\bibitem{FFTalgorithm}
F.~T. Luk and S.~Qiao, ``A fast eigenvalue algorithm for hankel matrices,''
  \emph{Linear Algebra and Its Applications}, vol. 316, no. 1-3, pp. 171--182,
  2000.

\end{thebibliography}

	\begin{appendix}
		
		\textit{Proof of Theorem \ref{Hankelvectorprod}}
		\begin{pf}
			Our proof is constructive. We first reduce the Hankel matrix-vector product into an equivalent Toeplitz matrix-vector product and eventually into a circulant matrix-vector product. We then show that the last product can be computed efficiently using FFT.
			\par First, we multiply the Hankel matrix $\mathbf{H}\in\mathbb{R}^{n \times m}$ by a matrix $\mathbf{\Pi}\in\mathbb{R}^{m\times m}$ to obtain a Toeplitz matrix $\mathbf{T}\in\mathbb{R}^{n\times m}$. The matrix $\mathbf{\Pi}$ is required to have the following structure:
			\begin{align}
			\mathbf{\Pi}=
			\begin{pmatrix}
			0 & 0 & \cdots &  0 & 1\\
			0 & 0 & \cdots &  1	& 0\\
			\vdots &  &\ddots & & \vdots \\
			1 & 0 & \cdots & 0 & 0 \\	
			\end{pmatrix}.
			\end{align}	
			One can easily verify that indeed:
			\begin{align}
			\mathbf{H}\cdot \mathbf{\Pi}=\mathbf{T} \label{HP=T}.
			\end{align}
			Using \eqref{HP=T}, we can express the product of the Hankel matrix $\mathbf{H}$ with the vector $v$ as:
			\begin{align}
			\mathbf{H}\cdot v= \mathbf{T}\cdot \mathbf{\Pi}^{-1} \cdot v.
			\end{align}
			One can additionally verify that:
			\begin{align}
			\mathbf{T}\cdot \mathbf{\Pi}=\mathbf{H}.
			\end{align}
			Thus, we also have that:
			\begin{align}
			\mathbf{H}\cdot v= \mathbf{T}\cdot\mathbf{\Pi} \cdot v=\mathbf{T}\cdot v_p.
			\end{align}
			where $v_p=\mathbf{\Pi}\cdot v$, $v_p\in\mathbb{R}^m$. This vector has the same elements as $v$ but sorted in reverse order:
			\begin{align}
			v_p&:=\begin{pmatrix} v_m & v_{m-1} & \cdots & v_1  \end{pmatrix}^\top.
			\end{align}
			So far, we have shown that the product  $\mathbf{H}\cdot v$ is equivalent to the product $\mathbf{T}\cdot v_p$. The next step is to embed the Toeplitz matrix $\mathbf{T}$ into a larger circulant matrix $\mathbf{C}\in\mathbb{R}^{(n+m-1)\times (n+m-1) }$ whose product with a vector can computed efficiently. We construct the matrix $\mathbf{C}$ as follows:
			\begin{align}
			\mathbf{C}=\begin{pmatrix}
			T & \star\\
			\star & \star
			\end{pmatrix}.\label{circulant}
			\end{align}
			We emphasize here that, a $n\times m$ Toeplitz matrix $\mathbf{T}$ should be embedded  in a $(n+m-1)\times(n+m-1)$ circulant matrix $\mathbf{C}$ with the matrix $\mathbf{T}$ being on its upper left block. This is because the distinct elements of the $(n+m-1)$ diagonals of the Toeplitz matrix $\mathbf{T}$ are those that define the vector $c$, which characterizes the circulant matrix $\mathbf{C}$ and precisely matches its first column.
			Moving forward, in light of \eqref{circulant}, the product $\mathbf{T}\cdot v_p$ can be expressed as a function of the circulant matrix $\mathbf{C}$ as follows:
			\begin{align}
			\mathbf{T} \cdot v_p=\begin{pmatrix} I_n & 0_{m-1} \end{pmatrix} \cdot \mathbf{C} \cdot v_e \label{toeplitzcirc},
			\end{align}
			where the vector $v_e$ is defined as:
			\begin{align}
			v_e=\begin{pmatrix}
			v_p\\
			0_{n-1}
			\end{pmatrix}.
			\end{align}
			We know that the circulant matrix $\mathbf{C}$ has the nice property of being diagonalized by the FFT matrix $\mathbf{F}$. That is:
			\begin{align}
			\mathbf{C}=\mathbf{F}^{-1} \mathbf{\Lambda} \mathbf{F}.
			\end{align}
			Hence, the product $\mathbf{C} \cdot v_e$ can be computed as:
			\begin{align}
			\mathbf{C}\cdot v_e=\mathbf{F}^{-1} \mathbf{\Lambda} \mathbf{F}v_e,
			\end{align}
			which can be written as:
			\begin{align}
			\mathbf{C}\cdot v_e=\mathbf{F}^{-1} \mathbf{\Lambda} \hat{V_e},
			\end{align}
			where $\hat{V_e}$ is the discrete Fourier transform of $v_e$, and $\mathbf{\Lambda}=\text{diag}(\mathbf{F}c)=\text{diag}(\hat{C})$ is a diagonal matrix with the elements of the discrete fourier transform of the vector $c$, $\hat{C}$, on its diagonal. From \eqref{toeplitzcirc}, it is easy to conclude that $\mathbf{H}\cdot v$ can be obtained as the first $n$ elements of the Hadamard product of the vectors $\hat{C}$ and $\hat{V}_e$. Note that $c$ is the vector that fully specifies the circulant matrix $\mathbf{C}$ and corresponds to its first column. That completes the proof.
		\end{pf}
		Algorithm \ref{algorithm1} can be used to efficiently compute  the product $p_{\text{H}}$ of a non-square Hankel matrix $\mathbf{H}\in\mathbb{R}^{n\times m}$ with a vector $v\in\mathbb{R}^m$. In this analysis, we computed the complexity of the proposed algorithm and proved  that it carries out the appropriate computation by exploiting FFT.
		

	\end{appendix}
\end{document}